\documentstyle[amssymb,11pt]{article}
\input{xypic}
\newcommand{\ortho}{{\perp}}

\newcommand{\calL}{{\cal L}}
\newcommand{\calK}{{\cal K}}
\newcommand{\calN}{{\cal N}}

\newcommand{\oka}{{\cal O}}
\newcommand{\CC}{{\Bbb C}}
\newcommand{\PP}{{\Bbb P}}
   
\newcommand{\ZZ}{{\Bbb Z}}

\newcommand{\NX}{{\calN_{X/Y}^\ast}}

\newcommand{\Thm}[1]{{\bf #1. $\:$}}
\renewcommand{\Text}[1]{\quad \mbox{\rm #1} \:\:}
\newcommand{\proof}{{\em Proof. $\:$}}

\newcommand{\ebew}{\hfill$\Box$ \par}
\newcommand{\Pic}{{\mbox{\rm Pic}\,}}
\newcommand{\OX}{\Omega_X}
\newcommand{\OY}{\Omega_Y}

\newtheorem{thm}{Theorem}
\newtheorem{prop}[thm]{Proposition}
\newtheorem{lemma}[thm]{Lemma}
\newtheorem{defi}[thm]{Definition}
\newtheorem{rk}[thm]{Remark}

\newtheorem{cor}[thm]{Corollary}

\begin{document}
\title{A Barth--Lefschetz theorem for toric varieties
}
\author{J\"org Zintl}
\date{ \small{\em  Dedicated to the memory of Michael Schneider}}
\maketitle


By the well known theorem of Barth--Lefschetz for complex projective space $Y
 = \PP^r$ one has
\[ H^q(Y, X; \CC) = 0 \quad \mbox{for} \:\: q \le 2n - r +1  \]
for any $n$-dimensional submanifold $X \subseteq Y$. 
This was proved first by Barth \cite{B} generalizing a theorem of Lefschetz on 
hypersurfaces $X$.
Using Le Potier's vanishing theorem \cite{LP} a new proof was given in a joint paper with Schneider 
\cite{SZ}. That article contains also
 references to further proofs of the theorem. Using the same arguments as in 
\cite{SZ}, Debarre \cite{D}
was able to show that the theorem of Barth--Lefschetz holds for abelian 
varieties $Y$ as well, assuming that the normal 
bundle of $X$ in $Y$ is ample.\footnote{
In fact, this is a special case of a result on submanifolds of projective homogeneous manifolds which has been proved earlier by Sommese \cite{So} using different methods.
}
\par
The idea behind this paper was to study the more general situation of a submanifold $X$ in a smooth projective toric variety $Y$. Giving up $\PP^r$ as ambient space, which is very special indeed, one needs some additional assumptions on the submanifold $X$. One requirement is that the normal bundle of $X$ in $Y$ is  ample, which is always true if $Y = \PP^r$. Besides that, the submanifold should meet the boundary of $Y$ transversally,  which again can be achieved in $\PP^r$ using the homogeneity of projective space. 
\par
For an arbitrary such submanifold $X$ in $Y$ one cannot expect that the result of Barth--Lefschetz still holds. In theorem \ref{6} a precise toric characterization is given in which cases it does. The main theorem \ref{mainthm} settles the general case. It gives a formula for all Hodge numbers $h^{p,q}(X)$ for any such submanifold $X$ in the Barth--Lefschetz range $p+q \le 2n-r$.  
\par
\begin{lemma}\label{1}
Let $ 0 \rightarrow F \rightarrow E_0 \rightarrow \ldots \rightarrow E_k 
\rightarrow 0$ be an exact sequence of sheaves on a scheme $Z$. Let $q \ge 0$ 
and assume 
\[ H^{q-i}(Z, E_i) = H^{q-i-1}(Z, E_i) = 0 \quad \mbox{for} \: \: 0 \le i \le 
k-1 .\]
Then
\[ H^q(Z, F) \cong H^{q-k}(Z, E_k) .\]
\end{lemma}
\par

\proof 
By induction on $k$ and cutting the sequence into two pieces.
\ebew
\par
Let $Y$ denote a smooth projective toric variety of dimension $r>0$, given by a 
fan $\Delta$ in some lattice $N \cong \ZZ^r$. The dual lattice to $N$ will be 
called  $M$. Each cone $\sigma \in \Delta$ corresponds to precisely one orbit of 
the action of $T_N  = (\CC^\ast)^r$ on $Y$. The closure of this orbit is denoted by 
$V(\sigma)$. It is again a smooth projective toric subvariety of codimension equal 
to the dimension of $\sigma$. Put
\[ \Delta(i) := \{ \sigma \in \Delta : \mbox{dim}\, \sigma = i \} .\]
For each $0 \le i \le r$  the set $\Delta(i)$ is finite because $Y$ is compact. Since $Y$ is 
smooth, Ishida's complex of degree  $p$ 
\[ 0 \rightarrow \Omega_Y^p \rightarrow \calK^0(Y; p) \rightarrow \ldots 
\rightarrow \calK^p(Y; p) \rightarrow 0 \]
is exact for $p = 0, \ldots, r$, where
\[ \calK^j(Y; p) = \bigoplus_{\sigma \in \Delta(j)} \oka_{V(\sigma)} \otimes_\ZZ
 \bigwedge^{p-j} (M \cap \sigma^\ortho ) .\]
Note that $M \cap \sigma^\ortho$ is a free $\ZZ$-module of rank $r-j$. 
\par
Let $X$ be a submanifold of $Y$ of dimension $n$ with ample normal bundle 
$\calN_{X/Y}$. We assume $X$ to be in general position. This is made precise in the following definition.

\begin{defi}\em
We say that $X \subset Y$ is {\em $\Delta$-transversal} if for all $\sigma \in \Delta(k)$ with $1 \le k \le n$ the  intersection of $X$ with $V(\sigma)$ is either empty or transversal, and if it is empty for $k > n$. 
\end{defi}

From now on we will assume that $X$ is $\Delta$-transversal. In particular this implies that $X \cap V(\sigma)$ is smooth, if it is not empty, and
\[ \calN_{X/Y} | X \cap V(\sigma) \: \cong \: \calN_{X \cap V(\sigma) / 
V(\sigma)} \]
is ample 
with rk $ \calN_{X/Y} | X \cap V(\sigma) = $\, rk $\calN_{X \cap V(\sigma) / 
V(\sigma)}$. 

We will need the following vanishing theorem due to Le Potier.

\begin{thm} {\rm (Le Potier)}
Let $Z$ be a smooth projective manifold of dimension $m$ and $E$ an ample vector 
bundle of rank $s$ on $Z$. Then 
\[ H^q(Z, \Omega_Z^p \otimes E^\ast ) = 0 \quad \mbox{for} \:\:  p+q \le m-s .\]
\end{thm}

In particular by what we said above, for all $\sigma \in \Delta(i)$ one has  
\[ H^q(X \cap V(\sigma), \calN_{X \cap V(\sigma) / V(\sigma)}^\ast ) = 0 \quad 
\mbox{for} \: q \le 2n-r-i . \]

We will use this to prove a vanishing theorem for symmetric powers of the 
conormal bundle $\NX$.

\begin{thm}\label{2}
Let $X$ be a $\Delta$-transversal submanifold of $Y$ with ample normal bundle $\calN_{X/Y}$. Then for all $k \ge 1$ 
\[ H^q(X, S^k \NX ) = 0 \Text{for}  q \le 2n-r .\]
\end{thm}
  
Note that if $X$ is a hypersurface in $Y$ then the theorem is true even without any condition on the position of $X$. This follows immediately from Kodaira's vanishing theorem.

\proof
Since $S^{k+1} \NX$ is a direct summand of $S^k \NX \otimes \NX$ it suffices to 
show
\[ H^q(X, S^k\NX \otimes \NX ) = 0 \Text{for}k \ge 0 \Text{and}  q \le 2n-r .\]
If $k=0$ this is just Le Potier's theorem. We proceed by induction on $k$. Take 
symmetric powers of the conormal seqence
\[ 0 \rightarrow \NX \rightarrow \OY^1|X \rightarrow \OX^1 \rightarrow 0 \]
to get
\[ 0 \rightarrow S^k\NX  \rightarrow S^{k-1}\NX\otimes\OY^1|X \rightarrow \ldots 
\rightarrow \OY^k|X \rightarrow \OX^k \rightarrow 0 . \]
Tensoring this with $\NX$ gives
\[ 0 \rightarrow S^k \NX \otimes \NX \rightarrow E_1 \rightarrow \ldots 
\rightarrow E_k \rightarrow \OX^k \otimes \NX \rightarrow 0 \]
where
\[ E_i = S^{k-i}\NX \otimes \OY^i|X \otimes \NX \Text{for} 1 \le i \le k. \]
It suffices to prove 
\[ H^{q-i+1}(X, E_i) = H^{q-i}(X, E_i) = 0 \Text{for} 1 \le i \le k \]
and $q \le 2n-r$, because then by lemma \ref{1} 
\[ H^q(X, S^k\NX \otimes \NX) \cong H^{q-k}(X, \OX^k \otimes \NX) \]
and we are reduced to Le Potier's theorem. 

So fix an integer $ 1 \le i \le k$. We may assume $i \le r$ since otherwise $E_i = 0$ and we are done. 
By tensoring Ishida's complex of degree  $i$  with $S^{k-i}\NX \otimes \NX$ we 
obtain
\[ \begin{array}{c} 
0 \rightarrow E_i \rightarrow \calK^0(Y; i) \otimes S^{k-i}\NX\otimes\NX 
\rightarrow \ldots \\
\qquad \ldots \rightarrow \calK^i(Y; i) \otimes S^{k-i}\NX\otimes\NX\rightarrow 
0 .
\end{array} \]
If $0 \le \nu \le i-1$ and $X \cap V(\sigma) \neq \emptyset$ for $\sigma \in \Delta(\nu)$, then  
$\NX | X \cap V(\sigma) \cong \calN_{{X\cap V(\sigma) / V(\sigma)}}^\ast$ 
as mentioned above, and one has 
\[ H^q(X \cap V(\sigma), S^{k-i} \calN_{{X\cap V(\sigma) / V(\sigma)}}^\ast 
\otimes  \calN_{{X\cap V(\sigma) / V(\sigma)}}^\ast ) = 0 \]
for $ q \le 2n-r-\nu$ by induction on $k$. Hence
\[ H^q(X, \calK^\nu(Y; i) \otimes S^{k-i}\NX \otimes \NX) = 0 \]
for $0 \le \nu \le i-1$ and $q \le 2n-r -\nu$.  This implies by lemma \ref{1}
\[ H^q(X, E_i) \cong H^{q-i}(X, \calK^i(Y; i) \otimes S^{k-i}\NX \otimes \NX)
 \]
for $q \le 2n-r$. Again by induction on $k$ this vanishes as we have $q-i \le 2n-r-i$. 
\ebew

As in \cite{SZ} we will compare the cohomology groups  $H^k(Y, \CC)$ and $H^k(X, \CC)$ 
using Hodge decomposition and factoring the natural map via
\[ H^i(Y,\OY^j) \rightarrow H^i(X, \OY^j|X) \rightarrow H^i(X, \OX^j) .\]

\begin{prop}\label{4}
Let $X \subset Y$ be as above. The natural map
\[  H^i(X, \OY^j|X) \rightarrow H^i(X, \OX^j) \]
is surjective if $i+j \le 2n-r$ and injective if $i+j \le 2n-r+1$.
\end{prop}

\proof
Consider the $j$-th symmetric power of the conormal sequence
\[ \begin{array}{rcl}
0 \rightarrow S^j \NX \rightarrow \ldots \rightarrow \NX\otimes \OY^{j-1}|X& 
\rightarrow &\OY^j|X \rightarrow \OX^j \rightarrow 0 \\
\searrow&&\nearrow\\
&M_j^\ast&\\
\nearrow&&\searrow\\
0 \qquad && \qquad 0
\end{array} \]
and cut it into two pieces. 
We have to prove $H^i(M_j^\ast) = 0$ for $i+j \le 2n-r+1$. Look at Ishida's 
complex of degree  $\nu$ tensored with $ S^{j-\nu} \NX $ for $1 \le \nu \le j-1$.
\[ \begin{array}{c}
 0 \rightarrow \OY^\nu|X  \otimes S^{j-\nu}\NX \rightarrow \calK^0(Y; \nu)\otimes 
S^{j-\nu}\NX \rightarrow \ldots \\[10pt]
\qquad \ldots \rightarrow \calK^\nu(Y; \nu)\otimes S^{j-\nu}\NX \rightarrow 0 .
\end{array}  \]
As in the proof of theorem \ref{2} we use the isomorphism  $\NX | X \cap 
V(\sigma) \cong \calN_{{X\cap V(\sigma) / V(\sigma)}}^\ast$ for $\sigma \in 
\Delta(\mu)$ if $X \cap V(\sigma) \neq \emptyset$. We obtain
\[ H^q(X, \calK^\mu(Y; \nu) \otimes S^{j-\nu} \NX) = 0 \Text{for} q \le 2n-r-\mu 
\]
and $0 \le \mu \le \nu-1$ by theorem \ref{2}. By the same theorem and by lemma 
\ref{1} we have
\[ H^q(  \OY^\nu|X \otimes  S^{j-\nu}\NX  ) \cong H^{q-\nu}(X, \calK^\nu(Y; \nu) 
\otimes S^{j-\nu} \NX) = 0 \]
if $q \le 2n-r$. 
Now we apply this and lemma \ref{1} to the long exact sequence containing 
$M_j^\ast$ above. We get
\[ H^{i+j-1}(X, S^j \NX) \cong H^i(X, M_j^\ast) \]
for $i+j-1 \le 2n-r$. Since the left hand side vanishes by theorem \ref{2}, we finally obtain $ H^i(X, M_j^\ast) = 0$ for $i+j \le 2n-r+1$ as desired. 
\ebew

Before we can proceed we need to show the following general fact.

\begin{prop}{\label{WZ}}
Let $Z$ be a smooth projective variety and let $W \subset Z$ be a submanifold with ample normal bundle $\calN_{W/Z}$. Let $D_1, \ldots, D_s$ be hypersurfaces on $Z$ such that $D_i \cap W = \emptyset$ for all $i=1, \ldots,s$. Assume that there are integers $a_1,\ldots, a_s \in \ZZ$ such that $\sum\limits_{i=1}^s a_i D_i$ is linearly equivalent to $0$. Then $a_i = 0$ for all $i = 1, \ldots,s$.
\end{prop}

\proof
Assume that the claim does not hold. We denote by $D_0$ the effective divisor $\sum a_i D_i$ where the sum runs over all $a_i > 0$, and by $D_\infty$ the effective divisor $\sum |a_i| D_i$ where $a_i < 0$. Since by assumption $D_0 - D_\infty$ is linearly equivalent to $0$, there exists a global meromorphic function $g$ on $Z$, such that
\[ (g) = D_0 - D_\infty .\]
Since $g$ has neither poles nor zeroes when restricted to $W,$ its restriction  has to be equal to a constant $c$, say. Hence the function $g-c$ vanishes on $W$, and it is defined
near $W$. 
Since it is not identically $0$, there is some integer
$k$ such that $g-c$ belongs
to ${\cal I}_W^k$ but $g-c$ does not belong
to
${\cal I}_W^{k+1}$, where ${\cal I}_W$ denotes the ideal sheaf of $W$ in $Z$. Hence $g-c$ provides a non-zero section of
${\cal I}_W^k/{\cal I}_W^{k+1}=S^k{\cal N}^\ast_{W/Z}$, which contradicts
the ampleness of ${\cal N}_{W/Z}$.
\ebew

\begin{cor}\label{corWZ}
Let $X \subset Y$ be as before. Let $s$ be the number of $\sigma \in \Delta(1)$ such that $V(\sigma) \cap X = \emptyset$. Then there is a short exact sequence
\[ 0 \rightarrow \OY^1|X \rightarrow \bigoplus_{{\stackrel{\sigma\in\Delta(1)}{V(\sigma)\cap X \neq \emptyset} }} \oka_X( -V(\sigma)\cap X ) \rightarrow \oka_X^{\delta(1)-r-s} \rightarrow 0. \]
\end{cor}

\proof
By \cite{BC} there is an exact generalized Euler sequence
\[ 0 \rightarrow \OY^1 \rightarrow \bigoplus_{\sigma\in\Delta(1)} \oka_Y( -V(\sigma) ) \rightarrow \oka_Y^{\delta(1)-r} \rightarrow 0 .\]
By construction the second map factors as
\[ \bigoplus_{\sigma\in\Delta(1)} \oka_Y(-V(\sigma)) \rightarrow
\bigoplus_{\sigma\in\Delta(1)} \oka_Y \rightarrow
\bigoplus_{i=1}^{\delta(1)-r} \oka_Y  \]
where the first map is injective, and the second map is obtained from the standard exact sequence
\[ 0 \rightarrow M \rightarrow \bigoplus_{\sigma\in\Delta(1)} \ZZ \,  V(\sigma) \rightarrow \Pic(Y) \rightarrow 0 \]
by tensoring with $\oka_Y$. Here $\bigoplus \ZZ \, V(\sigma)$ denotes the free group generated by the boundary divisors $V(\sigma)$, and we use $\Pic(Y) \cong \ZZ^{\delta(1)-r} $.
Now proposition \ref{WZ} implies that the map 
\[ \bigoplus_{{\stackrel{\sigma\in\Delta(1)}{V(\sigma)\cap X =\emptyset} }} \ZZ \, V(\sigma)  \rightarrow \Pic(Y) \]
is injective. This remains true when we tensor with $\oka_X$. Hence the gereralized Euler sequence restricted to $X$ splits off $s$ copies of $\oka_X$, one for each $\sigma \in \Delta(1)$ with $V(\sigma) \cap X = \emptyset$, and thus reduces to a short exact sequence
\[ 0 \rightarrow \OY^1|X \rightarrow \bigoplus_{{\stackrel{\sigma\in\Delta(1)}{V(\sigma)\cap X \neq \emptyset} }} \oka_X( -V(\sigma)\cap X ) \rightarrow \oka_X^{\delta(1)-r-s} \rightarrow 0 \]
as claimed.
\ebew

Using this observation, we can begin to determine Hodge numbers of $X \subset Y$. 

\begin{prop}\label{37}
Let $X \subset Y$ be as above and connected. Then for the Hodge numbers of $X$ holds
\[ h^{i,0}(X) = h^{0,i}(X) = 0 \]
for $1 \le i \le 2n-r$.
\end{prop}

\proof
Note that for $i \ge 1$ we have a natural inclusion
\[ \OY^i \hookrightarrow \OY^{i-1} \otimes \OY^1 \]
as a direct summand. Tensoring the exact sequence of corollary \ref{corWZ} with $\OY^{i-1}|X$ we get an injective map 
\[ \OY^{i-1}|X \otimes \OY^1|X \hookrightarrow \bigoplus_{{\stackrel{\sigma\in\Delta(1)}{V(\sigma)\cap X \neq \emptyset} }} \OY^{i-1}|X \otimes \oka_X(-V(\sigma)\cap X)  \]
Hence for each $\tau \in \Delta(1)$ there is a chain of injections
\[ H^0(X,\OY^{i}|X \otimes \oka_X(-V(\tau)\cap X)) 
\rightarrow H^0(X, \OY^{i}|X)
\rightarrow \qquad \qquad  \qquad \qquad \]
\[  \rightarrow  H^0(X,\OY^{i-1}|X \otimes \OY^1|X )
\rightarrow  \bigoplus_{{\stackrel{\sigma\in\Delta(1)}{V(\sigma)\cap X \neq \emptyset} }} H^0(X,\OY^{i-1}|X \otimes \oka_X(-V(\sigma)\cap X)) \]
If $i=1$ then the sum of cohomology groups on the right hand side vanishes, and so does the whole chain. Inductively this shows
\[ H^0(X,\OY^i|X) = 0 \]
for all $i \ge 1$. Because of proposition \ref{4} and Hodge symmetry the claim follows.
\ebew

\begin{defi}\em
Let $\delta(i)$ denote the number of elements in $\Delta(i)$. We write 
$\delta_{X,\sigma}$ for the number of connected components of $X\cap V(\sigma)$ 
for $\sigma \in \Delta$. Put
\[ \delta_X(i) := \sum_{\sigma \in \Delta(i)} \delta_{X,\sigma} .\]
\end{defi}

Using this notation we can state our main theorem.

\begin{thm}\label{mainthm}
Let $Y$ be a smooth projective toric variety of dimension $r$ and let $X$ be a $n$-dimensional, connected and $\Delta$-transversal submanifold of $Y$ with ample normal bundle. Then for $i+j \le 2n-r$ the following formula holds for the Hodge numbers of $X$. 
\[ h^i(X, \OX^j) = \left\{ \begin{array}{ll}
\sum_{k=0}^j (-1)^{j-k} {r-k \choose j-k} \delta_X(k) &\Text{if} i=j\\
0 & \Text{if} i \neq j .
\end{array} \right. \]
\end{thm}

\proof
If $i=j=0$ this is clear. If either $i=0$ or $j=0$ vanishing of the cohomology was proved in proposition \ref{37}. By proposition \ref{4} it is sufficient to determine the dimensions of $H^i(X,\OY^j|X)$ for $i+j \le 2n-r$. Let us assume $i > j$ first. Since proposition \ref{37} implies
 $H^q(X, \calK^\nu(Y; j) |X ) = 0$ for $1 \le q \le 
2n-r-\nu$ we have by lemma \ref{1} applied to Ishida's complex of degree $j$ and   
restricted to $X$
\[ H^i(X, \OY^j|X) \cong H^{i-j}(X, \calK^j(Y; j) | X ) = 0. \]
Because of proposition \ref{4} and Hodge symmetry we have
\[ H^i(X,\OY^j|X) \cong H^i(X, \OX^j) \cong H^j(X, \OX^i) \cong H^j(X, \OY^i|X) 
\]
for $i+j \le 2n-r$.  
Hence the proposition is true for $i < j$ as well.

It remains to look at the case $i=j$ and $i \neq 0$. In particular we have then $i \le \frac{2n-r}{2}$. If $i = 1$ the claim follows from the short exact sequence
\[0 \rightarrow \OY^1 | X \rightarrow \calK^0(Y; 1)|X \rightarrow \calK^1(Y; 1)|X 
\rightarrow 0 \]
and proposition \ref{37}. Otherwise consider Ishida's 
complex of degree $i\ge 2$ cut into pieces.
\[ \begin{array}{c}
0 \rightarrow \OY^i  \rightarrow \calK^0(Y; i) \rightarrow \calL^1 
\rightarrow 0 \\
0 \rightarrow \calL^1 \rightarrow \calK^1(Y; i) \rightarrow \calL^2 
\rightarrow 0 \\
\vdots\\
0 \rightarrow \calL^{i-1} \rightarrow \calK^{i-1}(Y; i) \rightarrow \calK^i(Y; 
i) \rightarrow 0 
\end{array} \]
The $\calL^\nu$ for $1 \le \nu \le i-1$ are inductively defined by this short exact sequences. We restrict them to $X$. 
Since $H^\mu(X, \calK^\nu(Y; i)|X) = 0$ for $ 1 \le \mu \le 2n-r-\nu$ we have 
\[ H^q(X, \OY^i|X) \cong H^{q-1}(X, \calL^1|X) \cong \ldots \cong H^{q-i+1}(X, 
\calL^{i-1}|X) \]
if $2n-r \ge q \ge i$. This vanishes if $q \neq i$ as we have seen just before. In other words
\[ H^q(X, \calL^p|X) = 0 \Text{if} 1 \le p \le i-1, \:\:  1 \le q \le 2n-r-p  \Text{and} p+q \neq i . \]
If $q=i$ we get 
\[ \begin{array}{rcl}
H^i(X, \OY^i|X) &\cong&H^{i-1}(X, \calL^1|X) \cong \ldots \cong H^1(X,\calL^{i-1}|X). 
\end{array} \]
We obtain also an exact sequence
\[ \begin{array}{c}
0 \rightarrow H^0(X,\calL^{i-1}|X) \rightarrow H^0(X, \calK^{i-1}(Y; i)|X) \rightarrow  \qquad \qquad \qquad \\[5pt]
\qquad \qquad \qquad \rightarrow H^0(X, \calK^{i}(Y; i)|X) \rightarrow H^1(X,\calL^{i-1}|X) \rightarrow 0 
\end{array} \] 
as well as
\[ \begin{array}{c}
0 \rightarrow H^0(X, \calL^{i-2}|X) \rightarrow H^0(X, \calK^{i-2}(Y; i)|X) \rightarrow H^0(X, \calL^{i-1}|X) \rightarrow 0 \\
 \vdots \\[5pt]
0 \rightarrow H^0(X, \calL^{1}|X) \rightarrow H^0(X, \calK^{1}(Y; i)|X) \rightarrow H^0(X, \calL^{2}|X) \rightarrow 0 \\[5pt]
0 \rightarrow H^0(X, \OY^i|X) \rightarrow H^0(X, \calK^{0}(Y; i)|X) \rightarrow H^0(X, \calL^{1}|X) \rightarrow 0 
\end{array} \]
By proposition \ref{37} the last exact sequence simply says
\[  H^0(X, \calK^{0}(Y; i)|X) \cong  H^0(X, \calL^{1}|X) .\]
Taking everything together we get a long exact sequence
\[ 0 \rightarrow  H^0(X, \calK^{0}(Y; i)|X) \rightarrow 
 H^0(X, \calK^{1}(Y; i)|X) \rightarrow \ldots \qquad \qquad \qquad \]
\[ \qquad \qquad \qquad 
\ldots \rightarrow 
 H^0(X, \calK^{i}(Y; i)|X) \rightarrow H^i(X, \OY^i|X) \rightarrow 0 .\]
So adding up all contributions finally gives
\[ \begin{array}{rcl}
h^i(X, \OY^i|X) &=& \sum_{\nu = 0}^i (-1)^{i-\nu} h^0(X, \calK^{\nu}(Y; i)|X)\\
&=& \sum_{\nu=0}^i (-1)^{i-\nu} {r - \nu \choose  i - \nu} \delta_X(\nu) .
\end{array} \]
\ebew

\begin{rk}\em
The vanishing statement of the theorem above can be proved without referring to the toric nature of $Y$. We only need to know  that $H^i(Y,\Omega_Y^j) = 0$ for $i \neq j$, and that the normal bundle of $X$ in $Y$ is ample. 
This fact, and the proof of the lemma below, were kindly pointed out
to me by Claire Voisin. 

Let $j : X \rightarrow Y$ denote the inclusion. Then there is a
commutative diagram
\[ \diagram
H^k(Y) \rrto^{[X]} \dto_{j^\ast} && H^{k+2r-2n}(Y) \dto^{j^\ast}\\
H^k(X) \rrto^{c_{r-n}} \urrto^{j_\ast} && H^{k+2r-2n}(X)
\enddiagram \]
where the horizontal maps are given by the cup-products with the
cohomology class $[X]\in H^{2r-2n}(Y)$ of $X$ in $Y$, and the top Chern class
$c_{r-n} := c_{r-n}(\calN_{X/Y})$ of the normal bundle of $X$,
respectively. This implies that $j_\ast$ is injective if the cup-product map $c_{r-n} \cup $ is injective.  
 Since the map $j_\ast$ repects the  Hodge decompositions on
$X$ and $Y$, and since $H^i(Y,\OY^j)=0$ for all $0\le i,j \le r$ with
$i\neq j$, the 
vanishing $H^i(X,\Omega_X^j)=0$ for $i \neq j$, and $ j \le 2n-r$  
 follows immediately from the lemma below. 
\end{rk}

\begin{lemma}
Let $j : W \rightarrow Z$ be the inclusion of an $n$-dimensional
submanifold $W$ into an $r$-dimensional smooth projective variety
$Z$, with ample normal bundle. 
Then the cup-product 
\[ \diagram
H^k(W) \rto^{c_{r-n}\quad} & H^{k+2r-2n}(W) 
\enddiagram \]
with the  top Chern class
$c_{r-n} := c_{r-n}(\calN_{W/Z})$ of the normal bundle of $W$ in $Z$ is
injective for $k \le 2n-r$.
\end{lemma}

\proof
Let $\calN_{W/Z}$ be the normal bundle of $W$ in $Z$. By assumption,
the tautological bundle $\oka_{\PP(\calN_{W/Z})}(1)$ on the projective bundle
$\pi: \PP(\calN_{W/Z}) \rightarrow X$ is ample. Hence we can apply the
hard Lefschetz 
theorem to it. If $h := c_1(\oka_{\PP(\calN_{W/Z})}(1))$, then we have that for $i
\le \dim(\PP(\calN_{W/Z})) -1 = r-2$ the cup-product map 
\[ \diagram H^i(\PP(\calN_{W/Z})) \rto^{h \: } &
H^{i+2}(\PP(\calN_{W/Z}))
\enddiagram \]
is injective. Furthermore we have a relation 
\[ 0 = \sum_{i=0}^{r-n} (-1)^{r-n-i}\, h^i.\,
\pi^\ast\left(c_{r-n-i}(\calN_{W/Z})\right) \]
in the cohomology of  $ \PP(\calN_{W/Z})$. Now let $\alpha \in H^k(X)$
and assume $c_{r-n}\cup \alpha = 0$. Lifting into the cohomology on
$\PP(\calN_{W/Z})$ gives the relation 
\[ \pi^\ast(c_{r-n}\cup \alpha) = \sum_{i=1}^{r-n} (-1)^{i-1}\, h^i.\,
\pi^\ast\left(c_{r-n-i}(\calN_{W/Z}) \cup \alpha\right) .\]
As mentioned above, if $k \le 2n-r$, then the cup-product with $h$ is injective on
$H^{k+2r-2n-2}(\PP(\calN_{W/Z}))$ by the Hard Lefschetz theorem. Therefore the
class $ \sum_{i=1}^{r-n} (-1)^{i-1}\, h^{i-1}.\,
\pi^\ast\left(c_{r-n-i}(\calN_{W/Z}) \cup \alpha \right)$ is equal to $0$ as
well. But since the expression of a cohomology class on
$\PP(\calN_{W/Z})$ by powers of $h$ in this form is
unique, we conclude that $\pi^\ast(\alpha) = 0$, and hence $\alpha = 0$. 
\ebew

Now let us return to the question of how the theorem of Barth--Lefschetz generalizes to smooth projective toric varieties. 

\begin{defi}\em
Let $X \subset Y$ be $\Delta$-transversal. 
We say that $X$ has the property (BL) if $\delta_{X,\sigma} =1$ for all $\sigma \in \Delta(t)$ with $2t \le 2n-r$ and $\delta_{X,\sigma} \ge 1$ if $2t = 2n-r+1$. For such an $X$ in particular $\delta_X(t) = \delta(t)$ for all $2t \le 2n - r$.
\end{defi}

\begin{rk}\em\label{11}
Note that if $Y = \PP^r$ we can always assume that  $X \subset Y$ of dimension $n$ has the property (BL). This is seen as follows. If $G = PGL(r)$ is the group of homogeneous
 transformations on $\PP^r$ then by Kleinman's theorem for any $V(\sigma)$ there is a nonempty open subset $U$ of $G$, such that every transformation of $X$ by an element of $U$ intersects $V(\sigma)$ transversally or not at all. Since there are only finitely many subvarieties $V(\sigma)$ in $Y$ there exists a homogeneous transformation which moves $X$ into $\Delta$-transversal position. Now 
for $\sigma \in \Delta(t)$ with  $2t \le 2n-r+1$ and $X \cap V(\sigma) \neq \emptyset$ put $V := V(\sigma)$. Since codim$_Y V = t$ we have 
codim$_Y X \cap V = r-n+t \le r$, hence $\delta_{X,\sigma} \ge 1$. Let $D_1$ and $D_2$ be two connected components of $X 
\cap V$. If $2t \le 2n-r$ then dim$\,D_1$ $+$ dim$\,D_2 = 2(n-t) \ge r$, so $X\cap V$ is connected. 
\end{rk}

\begin{rk}\em
It was conjectured by Hartshorne \cite{H} that each smooth $X \subset Y$ with ample normal bundle is connected if $2n - r \ge 0$. This would imply that  $\delta_{X,\sigma} \le 1$ if $\sigma \in \Delta(t)$ with $2t \le 2n-r$. If $X$ is a hypersurface this is always true. 
\end{rk}

\begin{thm}\label{6}
Let $Y$ be a smooth projective toric variety of dimension $r$ and let  $X$ be a 
$n$-dimensional submanifold with ample normal bundle which is $\Delta$-transversal. Then the natural map
\[ H^k(Y, \CC\,) \rightarrow H^k(X, \CC\,) \]
is surjective for  $k \le 2n- r$ and injective for $k \le 2n-r+1$ if and only if $X$ has the property (BL).
\end{thm}

\proof
We use the Hodge decompositions of $H^k(Y, \CC)$ and $H^k(X, \CC)$. Because of proposition \ref{4} it suffices to prove the theorem for the second map of the factorization
\[ H^i(Y, \OY^j) \rightarrow H^i(X, \OY^j|X) \rightarrow H^i(X, \OX^j) .\]
By theorem \ref{mainthm} it is enough to consider the case $i=j$, since $H^i(Y,\OY^j) = 0$ if $i \neq j$ and $H^i(X,\OY^j|X)= 0$ if $i \neq j$ and $i+j \le 2n-r$. 
Obviously $H^0(Y,\oka_Y) = H^0(X,\oka_X)=\CC$ if and only if $\delta_X(0)=1$. If $i \ge 1 $ we found in the proof of theorem \ref{mainthm} a long exact sequence that computed $H^i(X,\OY^i|X)$. Using the restriction map from the analogous exact sequence for the toric variety $Y$ to $X$  we get a commutative diagram
\[ \begin{array}{lcccccc}
0 \rightarrow &H^0(Y, \calK^{0}(Y; i)) & \rightarrow\ldots\rightarrow &H^0(Y, \calK^{i}(Y; i)) &\rightarrow&H^i(Y, \OY^i)& \rightarrow 0\\[7pt]
&\alpha_{i,0} \downarrow && \alpha_{i,i}   \downarrow && \rho \downarrow \\[7pt]
0\rightarrow &H^0(\calK^{0}(Y; i)|X) & \rightarrow\ldots\rightarrow& H^0( \calK^{i}(Y; i)|X) &\rightarrow& H^i( \OY^i|X)& \rightarrow 0
\end{array} \]
The property (BL) is equivalent to the property that the maps
\[ \alpha_{\nu,k} : H^0(Y, \calK^{k}(Y; \nu)) \rightarrow H^0(X, \calK^{k}(Y; \nu)|X) \]
are isomorphisms if $0 \le k \le \frac{2n-r}{2}$ and injective if $k = \frac{2n-r+1}{2}$. Note that $\alpha_{\nu,k}$ is injective or an isomorphism for some value of $\nu$ if and only if it is for all possible values of $\nu$. Using this, the theorem follows by an easy induction on $i$. 
\ebew

\begin{rk}\em 
$(i)$
Consider the special case $Y = \PP^r$. For any submanifold $X \subset Y$ the normal bundle is ample. So by what we said in remark \ref{11} this theorem implies the usual Barth--Lefschetz theorem for projective space.\\
$(ii)$
Because of Kleinman's theorem for homogeneous spaces this result is true for a   submanifold $X \subset Y$ with ample normal bundle without any condition on $\Delta$-transversality if $Y$ is equal to a product of projective spaces.
\end{rk}

\Thm{Acknowledgements}
I am grateful to the Graduiertenkolleg ``Komplexe Mannigfaltigkeiten'' at the University of Bayreuth as well as to the University of Bath for their hospitality. In particular I wish to thank Dr. G.K. Sankaran for many helpful conversations. Partially this
 work has been supported by the EU HCM project ``Algebraic Geometry in Europe'' (AGE), contract number ERBCHRXCT 940557.

\noindent
J\"org Zintl\\
Fachbereich Mathematik\\
Universit\"at Kaiserslautern\\
Postfach 3049\\
67653 Kaiserslautern\\
Germany\\
{\tt zintl@mathematik.uni-kl.de}


\begin{thebibliography}{BC}
\bibitem[B]{B}{Barth W., Transplanting cohomology classes in complex-projective space, {\em Amer. J. Math.} 92 (1970) 951--967}
\bibitem[BC]{BC}{Batyrev V., Cox D., On the Hodge structure of projective hypersurfaces in toric varieties, {\em Duke Math. J.} 75 No. 2 (1994) 293--338}
\bibitem[D]{D}{Debarre O., Fulton--Hansen and Barth--Lefschetz theorems for subvarieties of abelian varieties, {\em Crelle} 467 (1995) 187--197}
\bibitem[H]{H}{Hartshorne R., {\em Ample Subvarieties of Algebraic Varieties}, LNM 156, Springer, Berlin, Heidelberg, New York (1970)}
\bibitem[O]{O}{Oda T., {\em Convex Bodies and Algebraic Geometry}, Erg. Math. 
Grenzgeb. 3.Folge, Band 15, Springer, Berlin, Heidelberg, New York (1988)}
\bibitem[LP]{LP}{Le Potier J., Annulation de la cohomologie \`a valeurs dans un 
fibr\'e vectoriel holomorphe positif de rang quelconque, {\em Math. Ann.} 218 
(1975) 35--53}
\bibitem[So]{So}{Sommese A.J., Complex subspaces of homogeneous complex manifolds II - Homotopy results, {\em Nagoya Math. J.} 86 (1982) 101--129}
\bibitem[SZ]{SZ}{Schneider M., Zintl J., The theorem of Barth--Lefschetz as a 
consequence of Le Potier's vanishing theorem, {\em manuscripta math.} 80 (1993) 
259--263}
\end{thebibliography}
\end{document}